\documentclass[12pt]{amsart}
\usepackage{amssymb}
 \hoffset -1.2cm \voffset -1.3cm

\newtheorem{theorem}{Theorem}[section]
\newtheorem{lemma}[theorem]{Lemma}
\newtheorem{proposition}[theorem]{Proposition}
\newtheorem{corollary}[theorem]{Corollary}
\theoremstyle{definition}
\newtheorem{definition}[theorem]{Definition}

\theoremstyle{remark}
\newtheorem{remarks}[theorem]{Remarks}
\newtheorem{remark}[theorem]{Remark}
\numberwithin{equation}{subsection}
\newfont{\ajb}{eufm10 at12pt}
\newfont{\aj}{eufm10 at10pt}
\newfont{\ajk}{eufm10 at8pt}
\newfont{\kh}{msbm10 at10pt}
\newfont{\khk}{msbm10 at 8pt}
\newcommand{\bea}{\begin{eqnarray*}}
\newcommand{\eea}{\end{eqnarray*}}
\begin{document}
\title[
 Skew product dynamical systems] {Skew product dynamical systems, Ellis groups and topological centre}
\author{A. Jabbari}
\address{Department of Mathematics, Ferdowsi University, P. O. Box 1159, Mashhad 91775, Iran}
\email{shahzadeh@math.um.ac.ir} \subjclass[2000]{37B05 ,54H20,
22A20, 43A60}\keywords{ distal flow, Ellis group, compact Hausdorff
right topological group, topological centre}
\author{H. R. E.  Vishki}
\address{Department of Mathematics, Ferdowsi University, P. O. Box 1159, Mashhad 91775,
Iran; Centre of Excellence in Analysis on Algebraic Structures
(CEAAS), Ferdowsi University of Mashhad, Iran.}
\email{vishki@ferdowsi.um.ac.ir}

\thanks{}

\subjclass[2000]{37B05, 37A05, 43A60}

\keywords{flow, minimal flow, distal flow, uniquely ergodic flow,
Ellis group, compact right topological group, topological centre}

\date{}

\dedicatory{}

\commby{}

\begin{abstract}
In this paper, a general construction of a skew product dynamical
system, for which the skew product dynamical system studied by Hahn
is a special case, is given. Then the ergodic and topological
properties (of a special type) of our newly defined systems (called
Milnes type systems) are investigated. It is shown that the Milnes
type systems are actually natural extensions of dynamical systems
corresponding to some special distal functions. Finally, the
topological centre of the Ellis group of any skew product dynamical
system is calculated.
\end{abstract}

\maketitle

\maketitle
\section {Introduction}

In his well-known work~\cite{H1}, Hahn studied ergodic and
topological properties of dynamical systems $(X,T)$ in which $X=G^k$
and $T:X\rightarrow X$ is defined by
\[T(x_1,x_2,\ldots,x_k)=(\gamma_0+x_1,x_1+x_2,\ldots,x_{k-1}+x_k)\] where $k=1,2,3,\ldots$ or $k=\infty$ and $G$ is a compact
monothetic group with generator $\gamma_0$, that is, $G$ is a
compact group and the subgroup generated by $\gamma_0$ is dense in
$G$ (so $G$ is
 abelian).
He called $(X,T)$ a skew product dynamical system. Among other
things, Hahn~\cite{H1} showed that if the dual group of $G$ has no
element of finite order, then $(X,T)$ is a (distal) minimal
dynamical system which is also uniquely ergodic. Such
transformations on the torus were first studied by Anzai in
~\cite{A}, and the question of strict ergodicity of these
transformations on the torus was studied by
Furstenberg~\cite{F}.\\
Also Milnes~\cite[Example 3]{M} studied a different (kind of skew
product)
dynamical system on the 4-torus defined as\\
\centerline{$T(x,y,z,u)=(\lambda^{4}x,\lambda^{6}x^{3}y,
\lambda^{4}x^{3}y^{2}z,\lambda xyzu)$} where $\lambda$ is an
irrational (i.e. non-root of unity) element of the unit circle.\\

In this paper, a general construction of a skew product dynamical
system is given. It is shown that Hahn's type skew product dynamical
systems~\cite{H1} and also the Milnes dynamical system on the
4-torus~\cite{M} are special cases of our definition of a skew
product system, (see section 3 for the definition and results).
Then, in section 4, we examine the ergodic and topological
properties of Milnes type skew product dynamical systems. These skew
product systems, on the finite dimensional tori, are actually
extensions of the dynamical systems corresponding to the functions
$n\rightarrow\lambda^{n^k}$ on the integers (Theorem~\ref{2.4}).
Also section 5 is devoted to a generalization of a result of E.
Salehi~\cite{S} on the unique ergodicity of the functions
$n\rightarrow\lambda^{n^k}$ ($k\in{\mathbb{N}}$,
$\lambda\in{\mathbb{T}}$) on the additive group of integers by using
the unique ergodicity of Milnes type skew product
systems.\\

An interesting problem concerning any compact right topological
semigroup is characterizing its topological centre. Such semigroups
in topological dynamics arise naturally in two examples. First,
semigroup compactifications of any semitopological semigroup, and
second, the enveloping semigroups of flows. In either case, there
are rather few situations in which the corresponding topological
centre is calculated explicitly. We refer the interested reader
to~\cite{LP},~\cite{N} or, in general,~\cite{DL} for some
interesting results on the first case. For the latter case, there is
a nice result due to Namioka~\cite{N2} which computes the
topological centre of the Ellis group of a dynamical system on the
2-torus. Section 6 of this paper is devoted to the study of the
Ellis groups of Hahn's type and Milnes type skew product dynamical
systems and their topological centres. In particular, our results
generalize Namioka's result.
\section{Preliminaries}
For the background materials and notations we follow Berglund et
al.~\cite{BJM} as much as possible. For a semigroup $S$,  the right
translation $\rho_{t}$ and the left translation $\lambda_{s}$ on $S$
are defined by $\rho_{t}(s)=st=\lambda_{s}(t)$, ($s, t\in S$). A
semigroup $S$, equipped with a topology, is said to be right
topological if all of the right translations are continuous,
semitopological if all of the left and right translations are
continuous.  If $S$ is a right topological semigroup then the set
$\Lambda(S)=\{s\in S: \lambda_{s}$ is continuous$\}$ is called the
topological centre of $S$.

 By a flow $(S,X)$ we mean a
semitopological semigroup $S$ and a compact Hausdorff space $X$ in
which there exists a mapping $S\times X\rightarrow X, (s,x)\mapsto
sx$ such that for all $s, t\in S$ and $x\in X$, $(st)x=s(tx)$ and
$ex=x$ if $S$ has an identity $e$, and the mapping $\varepsilon_{s}
:X\rightarrow X$, $x\mapsto sx$ is continuous for each $s\in S$. If
$(S,X)$ is a flow, then the closure of the set
$\{\varepsilon_s:X\rightarrow X:s\in S\}$ in $X^X$, with the
relativization of the product topology from $X^X$, is a compact
right topological semigroup under the composition operator. This
closure is denoted by $\Sigma(S,X)$, or simply $\Sigma$ if there is
no risk of confusion, and is called the enveloping semigroup of the
flow. It is obvious from the definition that $\Lambda(\Sigma(S,X))$
always contains the set $\{\varepsilon_s:s\in S\}$. The idea of
enveloping semigroup was introduced by Ellis~\cite{E2}. It was also
proved by Ellis~\cite{E1} that the enveloping semigroup of a flow
$(S,X)$ is a group if and only if $(S,X)$ is distal, that is,
$\lim_{\alpha} s_{\alpha}x =\lim_{\alpha} s_{\alpha}y$ for some net
$\{s_{\alpha}\}$ in $S$ and for some $x, y \in X$ implies that
$x=y$. The enveloping semigroup of a distal flow is called the Ellis
group of the flow. A flow $(S,X)$ is called minimal if the orbit
closure of each point is dense in $X$.

 The next proposition is a crucial tool
in the characterization of the topological centre of the Ellis group
corresponding to any minimal flow, and is used in section 6 to
characterize the topological centre of the Ellis groups of skew
product dynamical systems.
\begin{proposition}
\label{m} Let $(S,X)$ be a minimal flow, then the set ${\Sigma}_{c}$
of all continuous elements of $\Sigma$ coincides with
$\Lambda(\Sigma)$.
\end{proposition}
\begin{proof} Obviously ${\Sigma}_{c}\subseteq {\Lambda}(\Sigma)$.
Let ${\sigma}\in{\Lambda}(\Sigma)$, $x\in X$ and let
$\{x_{\alpha}\}$ be a net in $X$ such that $x_{\alpha}\rightarrow
x$. We have to show that ${\sigma}(x_{\alpha})\rightarrow
{\sigma}(x)$. Since $X$ is minimal $\Sigma x=X$, hence for each
$\alpha$ there exists $\tau_{\alpha}\in\Sigma$ such that
$\tau_{\alpha}x=x_{\alpha}$. By taking a convergent subnet of
$\tau_{\alpha}$ if necessary, we may assume that
$\tau_{\alpha}\rightarrow \tau\in\Sigma$; it follows that
$\tau_{\alpha}x \rightarrow \tau x$. On the other hand
$\tau_{\alpha}x=x_{\alpha}\rightarrow x$, hence $\tau x=x$. Since
${\sigma}\in{\Lambda}(\Sigma)$, $\sigma
x_{\alpha}=\sigma(\tau_{\alpha}x)=\sigma\circ\tau_{\alpha}(x)\rightarrow
\sigma\circ\tau(x)=\sigma(\tau x)=\sigma x$, as required.
\end{proof}
A dynamical system is a pair $(X,T)$ where $X$ is a compact
Hausdorff space and $T$ is a homeomorphism from $X$ onto $X$.
Obviously, if $(X,T)$ is a dynamical system, then $({\mathbb{Z}},X)$
is a flow with the action $(n,x)\mapsto T^n(x)$, for
$n\in{\mathbb{Z}}$ and $x\in X$. Hence $(X,T)$ is called minimal
(distal) if the corresponding flow $({\mathbb{Z}},X)$ is minimal
(distal).

 A dynamical system $(X,T)$ is called uniquely ergodic if
there exists a unique $T$-invariant positive probability measure
$\mu$ on $X$. A minimal uniquely ergodic dynamical system is called
strictly ergodic. A sequence $\{x_n\}$ in $X$ is said to be
uniformly distributed if for all $x\in X$ and for each $f\in C(X)$,
$lim_{N\rightarrow\infty}\frac{1}{N}\sum_{n=0}^{N-1}f(T^n(x))=\int_X
f d\mu$.

 The following two theorems due to Weyl~\cite{W} and
Oxtoby~\cite{O} will be used in section 4 to prove the minimality
and unique ergodicity of the Milnes type skew product dynamical
systems.
\begin{theorem}[Weyl]\label{Weyl}  Let
$p(z)=a_0+a_{1}z+\cdots+a_{k-1}z^{k-1}+a_{k}z^k$ be a polynomial of
degree $k$ with real coefficients. If for some $j$, $1\leq j\leq k$,
$a_j$ is irrational, then:
\[\sum_{n=0}^{N}\exp(p(n))=\circ(N)\]
 uniformly in
$a_0,a_1,\ldots,a_{k-1}$.
\end{theorem}

\noindent For an interesting proof of the following theorem
see~\cite{H1}.
\begin{theorem}[Oxtoby]\label{Oxtoby} Let $(X,T)$ be a minimal
dynamical system. If for each $f\in C(X)$
\[\lim_{N\rightarrow\infty}\frac{1}{N}\sum_{n=0}^{N-1}f(T^n(x))\]
exists uniformly in $x$, then $(X,T)$ is strictly ergodic.
\end{theorem}
\section{The construction of skew product dynamical systems}
If $\{J_i:i\geq 0\}$ is a sequence of non-zero real numbers with
$J_0=1$, define $J(i,0)=1$ and $J(i,1)=J_{i-1}$ for $i\geq 1$ and
\[J(i,l)=J_{i-1}J_{i-2}\cdots J_{i-l}/(J_1J_2\cdots
J_{l-1})~for~ i\geq l\geq 2.\] Let $Q_0=1$ and for each $i\geq 1$,
let $Q_i(n)$ be a polynomial which satisfies the following equation:
\begin{equation} \label{e:0}\sum_{j=l}^{i}J(i-l,j-l)Q_{i-j}(1)Q_{j-l}(n-1)=Q_{i-l}(n)~for~all~
i\geq l\geq 0~and~ for~all~ n\in{\mathbb{N}}\end{equation} For
example if for each $i$, $J_i=1$, then it is enough to take
$Q_i(n)=\frac{1}{i!}n(n-1)\cdots(n-i+1)$ for $i\leq n$ and
$Q_i(n)=0$ elsewhere. Also if for each $i$, $J_i=i+1$, then
$Q_i(n)=n^i$ satisfies the cited equation. Now we are ready to make
our main construction:
\begin{definition} \label{1.1}
Let $\{J_i:i\geq 0\}$ be a sequence of non-zero real numbers with
$J_0=1$ and let $Q_i(n)$ be polynomials satisfying \ref{e:0} such
that for $0\leq j\leq i\leq k$ one has
$J(k-j,i-j)Q_{i-j}(1)\in{\mathbb{Z}}$. Let $k\in{\mathbb{N}}$, let
$G$ be a compact monothetic group with generator $\gamma_0$. Put
$X=G^k$ and define $T:X\rightarrow X$ by
\[T(x_1,x_2,\ldots,x_k)=(y_1,y_2,\ldots,y_k)\] where for each
$i=1,2,\ldots,k$ \begin{equation} \label{e:01}
y_i=J(k,i)Q_i(1)\gamma_0+\sum_{j=1}^{i}J(k-j,i-j)Q_{i-j}(1)x_j
\end{equation} Then $(X,T)$ is called a skew product dynamical system.
\end{definition}
\noindent{\bf Two types of skew product systems.} (i) For each
$i\geq 0$ let $J_i=1$. Let $Q_i(n)=\frac{1}{i!}n(n-1)\cdots(n-i+1)$
for $1\leq i\leq n$ and $Q_i(n)=0$ for $i>n$. Then it is easily seen
that (in~\ref{e:01}) $y_1=\gamma_0+x_1$ and $y_i=x_{i-1}+x_i$ for
$i\geq 2$. Hence we get the Hahn's skew product transformation:
\[T(x_1,x_2,\ldots,x_k)=(\gamma_0+x_1,x_1+x_2,\ldots,x_{k-1}+x_k).\]
In this paper we refer to such systems as Hahn's type skew product
systems. In this case $k=\infty$ is also possible.\\
(ii) For each $i\geq 0$ let $J_i=i+1$ and let $Q_i(n)=n^i$. Some
computation reveals that for each $l\leq i$, $J(i,l)=(^{i}_{l})$ and
so if
$T(x_1,x_2,\ldots,x_k)=(y_1,y_2,\ldots,y_k)$, then\\
\centerline {$y_1=k\gamma_0+x_1$}\\
\centerline {$y_2=(^{k}_{2})\gamma_0+
(^{k-1}_{1})x_1+x_2$}\\
\centerline {$\vdots$}\\
\centerline {$y_i=(^{k}_{i})\gamma_0+\sum_{j=1}^{i}(^{k-j}_{i-j})x_j$}\\
\centerline {$\vdots$}\\
\centerline {$y_k=\gamma_0+\sum_{j=1}^{k}x_j$.}\\
In this case we call $(X,T)$ a Milnes type skew product
transformation, just because for $k=4$ and $G={\mathbb{T}}$ one has
$T(x_1,x_2,x_3,x_4)=(\lambda^{4}x_1,\lambda^{6}x_1^{3}x_2,
\lambda^{4}x_1^{3}x_2^{2}x_3,\lambda x_1x_2x_3x_4)$ which is nothing
but the Milnes dynamical system.
\section{Milnes type skew product dynamical systems}
The results and methods of this section are based on those of
Hahn~\cite{H1}. Unless otherwise stated, $(X,T)$ denotes a Milnes
type skew product dynamical system with the base group $G$, where
$G$ is a compact monothetic group with generator $\gamma_0$. The
following theorem is the main result of this section.
\begin{theorem}\label{1.7}
If the dual group of $G$ has no element of finite order, then

 (i) For each $x$
in $X$ the points $T^n(x)$, ($n=1,2,3,\ldots$) are uniformly
distributed in $X$, in particular $(X,T)$ is minimal.

(ii) $(X,T)$ is strictly ergodic. \end{theorem}

To prove our main Theorem we need to mention some preliminaries.
\begin{lemma} \label{1.2} Let $n\in{\mathbb{Z}}$. Then for
$x=(x_1,x_2,\ldots,x_k)$ in $X$,
$T^n(x_1,x_2,\ldots,x_k)=(y_{1,n},y_{2,n},\ldots,y_{k,n})$, where
for each $i=1, 2,\ldots, k$,
$y_{i,n}=(^{k}_{i})n^i\gamma_0+\sum_{j=1}^{i}(^{k-j}_{i-j})n^{i-j}x_j$.
\end{lemma}
\begin{proof} We prove the result for positive integers by induction,
the general result will be then observed easily. Clearly the result
holds for $n=1$. Let $n\neq 1$ be a natural number and assume that
for each $x=(x_1,x_2,\ldots,x_k)$ in $X$,
$T^{n-1}(x)=(y_{1,n-1},y_{2,n-1},\ldots,y_{k,n-1})$, in which
$y_{i,n-1}=(^{k}_{i}){(n-1)}^i\gamma_0+\sum_{j=1}^{i}(^{k-j}_{i-j}){(n-1)}^{i-j}x_j$
for each $i=1, 2,\ldots, k$. Let $x=(x_1,x_2,\ldots,x_k)$, so
$T^n(x)=T(T^{n-1}(x))=T
(y_{1,n-1},y_{2,n-1},\ldots,y_{k,n-1}):=(w_1,w_2,\ldots,w_k)$, hence
for each $i=1, 2,\ldots, k$ \bea
w_i&=&(^{k}_{i})\gamma_0+\sum_{j=1}^{i}(^{k-j}_{i-j})y_{j,n-1}\\
&=&(^{k}_{i})\gamma_0+\sum_{j=1}^{i}(^{k-j}_{i-j})[(^{k}_{j}){(n-1)}^j\gamma_0+\sum_{l=1}^{j}(^{k-l}_{j-l}){(n-1)}^{j-l}x_l]\\
&=&(^{k}_{i})\gamma_0+\sum_{j=1}^{i}(^{k-j}_{i-j})(^{k}_{j}){(n-1)}^j\gamma_0+\sum_{j=1}^{i}
\sum_{l=1}^{j}(^{k-j}_{i-j})(^{k-l}_{j-l}){(n-1)}^{j-l}x_l\\
&=&[\sum_{j=0}^{i}(^{k-j}_{i-j})(^{k}_{j}){(n-1)}^j]\gamma_0+\sum_{l=1}^{i}\sum_{j=l}^{i}(^{k-j}_{i-j})(^{k-l}_{j-l}){(n-1)}^{j-l}x_l\\
&=&[\sum_{j=0}^{i}(^{k}_{i})(^{i}_{j}){(n-1)}^{j}]\gamma_0+{\sum_{l=1}^{i}\sum_{j=l}^{i}(^{k-l}_{i-l})(^{i-l}_{j-l}){(n-1)}^{j-l}}x_l,
 ((^{k-j}_{i-j})(^{k-l}_{j-l})=(^{k-l}_{i-l})(^{i-l}_{j-l}))\\
&=&(^{k}_{i})[\sum_{j=0}^{i}(^{i}_{j}){(n-1)}^{j}]\gamma_0+{\sum_{l=1}^{i}(^{k-l}_{i-l})[\sum_{j=l}^{i}(^{i-l}_{j-l}){(n-1)}^{j-l}}]x_l\\
&=&(^{k}_{i})n^i\gamma_0+\sum_{l=1}^{i}(^{k-l}_{i-l})n^{i-l}x_l,(
\sum_{j=l}^{i}(^{i-l}_{j-l}){(n-1)}^{j-l}=n^{i-l}). \eea The lemma
is now established.
\end{proof}

\begin{remark}
If $(X,T)$ is any skew product dynamical system, then with the same
method as in the above Lemma, and the notations of section 3, one
may show that: For $n\in{\mathbb{Z}}$ and  $x=(x_1,x_2,\ldots,x_k)$
in $X$, $T^n(x_1,x_2,\ldots,x_k)=(y_{1,n},y_{2,n},\ldots,y_{k,n})$,
where for each $i=1, 2,\ldots, k$,
$y_{i,n}=J(k,i)Q_i(n)\gamma_0+\sum_{j=1}^{i}J(k-j,i-j)Q_{i-j}(n)x_j$.
\end{remark}
\begin{corollary} \label{1.3} $(X,T)$ is distal.
\end{corollary}
\begin{proof} Let $(x_1,x_2,\ldots,x_k)$ and
$(x_1^{'},x_2^{'},\ldots,x_k^{'})$ be members of $X$ such that
\[\lim_\alpha T^{m_\alpha} (x_1,x_2,\ldots,x_k)=\lim_\alpha
T^{m_\alpha}(x_1^{'},x_2^{'},\ldots,x_k^{'})\] for some net
$\{m_{\alpha}\}$ in ${\mathbb{Z}}$. Then from Lemma ~\ref{1.2} for
each $i=1,2,\ldots,k$ we obtain
\[\lim_{\alpha}[(^{k}_{i}){m_{\alpha}^{i}}\gamma_0+\sum_{j=1}^{i}{{(^{k-j}_{i-j})}m_\alpha^{i-j}}x_{j}]=
\lim_{\alpha}[{(^{k}_{i}){m_{\alpha}^{i}}}\gamma_0+\sum_{j=1}^{i}{{(^{k-j}_{i-j})}m_\alpha^{i-j}}{x_{j}^{'}}]
.\]  Let $i=1$, then
$\lim_{\alpha}{km_{\alpha}}\gamma_0+x_1=\lim_{\alpha}{km_{\alpha}}\gamma_0+x_1^{'}$,
hence $x_1=x_1^{'}$. Now let $1\leq l<k$ and assume that we have
already shown that $x_j=x^{'}_{j}$ for $j=1,2,\ldots,l$. We will
show that $x_{l+1}=x^{'}_{l+1}$. Take $i=l+1$ in the above equation,
then we obtain
$\lim_{\alpha}{(^{k}_{l+1}){m_{\alpha}^{l+1}}}\gamma_0+x_{l+1}=\lim_{\alpha}{(^{k}_{l+1}){m_{\alpha}^{l+1}}}\gamma_0+x^{'}_{l+1}$.
Hence $x_{l+1}=x^{'}_{l+1}$. It follows by induction that
$(x_1,x_2,\ldots,x_k)=(x_1^{'},x_2^{'},\ldots,x_k^{'})$ that is,
$(X,T)$ is distal.
\end{proof}
 For
each $x$ in $X$, $T(x)=\gamma+U(x)$ where
$\gamma=(k\gamma_0,(^{k}_{2})\gamma_0,(^{k}_{3})\gamma_0\ldots,\gamma_0)$
and $U$ is the homomorphism on $X$ defined by $U=T-\gamma$. Hence
$T^n=((^{k}_{1})n\gamma_0,(^{k}_{2})n^2\gamma_0,\ldots,n^k\gamma_0)+U^n$,
for all $n\in{\mathbb{Z}}$. Let $\hat{G}$ and $\hat{X}$ denote the
dual groups of $G$ and $X$, respectively. And let $\hat{T}$ and
$\hat{U}$ denote the corresponding dual transformations of $T$ and
$U$, respectively.

\begin{lemma} \label{1.4} For each
$\eta=(\eta_1,\eta_2,\ldots,\eta_k)\in\hat{X}$ and
$n\in{\mathbb{Z}}$;

(i) $\hat{U}^n(\eta)=(\zeta_{1,n},\zeta_{2,n},\ldots,\zeta_{k,n})$
where $\zeta_{j,n}=\prod_{i=j}^{k}{\eta_i}^{(^{k-j}_{i-j})n^{i-j}}$
for $j=1,2,\dots,k$.

(ii)
$\hat{T}^n(\eta)=\eta((^{k}_{1})n\gamma_0,(^{k}_{2})n^2\gamma_0,\ldots,n^k\gamma_0)\hat{U}^n(\eta)$.
\end{lemma}
\begin{proof} (i)
Let $x\in X$ and let $U(x)=(z_1,z_2,\ldots,z_k)$, then
$z_i=\sum_{j=1}^{i}(^{k-j}_{i-j})x_j$. Now we have
$\hat{U}(\eta)(x)=\eta(Ux)=\eta_1(z_1)\eta_2(z_2)\cdots\eta_k(z_k)=\eta_1(x_1)
\eta_2((^{k-1}_{1})x_1+x_2)\cdots\eta_k(x_1+x_2+\ldots+x_k)=\zeta_1(x_1)\zeta_2(x_2)\cdots\zeta_k(x_k)$
where $\zeta_j=\prod_{i=j}^{k}{\eta_i}^{(^{k-j}_{i-j})}$ for
$j=1,2,\dots,k$. Hence the result holds for $n=1$. The general
result holds by induction. (ii) follows from (i) and the formulas
preceding the lemma.
\end{proof}

Let $\mu$ be the Haar measure on $X$. To prove our result on the
unique ergodicity of $(X,T)$, we first show that $(X,T)$ is ergodic
provided $\hat{G}$ has no (non-trivial) element of finite order:
\begin{theorem}\label{1.5} If the character group of $G$ has no
element of finite order, then $(X,\mu,T)$ is ergodic.
\end{theorem}
\begin{proof} Assume that $\hat{G}$ has no element of finite
order, and that $(X,\mu,T)$ is not ergodic. Hence by~\cite[Corollary
3]{H2} $\hat{U}^n\eta=\eta$ for some positive integer $n$ and for
some non-trivial character $\eta\in\hat{X}$. Let
$\eta=(\eta_1,\eta_2,\ldots,\eta_k)$ in which $\eta_i\in\hat{G}$,
then for each $x=(x_1,x_2,\ldots,x_k)$ in $X$,
$\eta(x)=\eta_1(x_1)\eta_2(x_2)\cdots\eta_k(x_k)$. By Lemma
~\ref{1.4} (i) the equation $\hat{U}^n\eta=\eta$ implies that
$\eta_{k-1}\eta_k^n=\eta_{k-1}$. Since $\hat{G}$ has no element of
finite order, $\eta_k\equiv 1$. Now by downward induction assume
that for some $i$ with $2<i\leq k$,
$\eta_i\equiv\eta_{i+1}\equiv\cdots\equiv\eta_k\equiv 1$. From
$\hat{U}^n\eta=\eta$ we conclude that
$\eta_{i-2}\eta_{i-1}^{(k-i+1)n}\equiv\eta_{i-2}$. Since $\hat{G}$
has no element of finite order we get $\eta_{i-1}\equiv 1$, so that
by induction $\eta_2\equiv\eta_3\equiv\cdots\equiv\eta_k\equiv 1$.
It remains to show that $\eta_1\equiv 1$. By ~\cite[Theorem 4]{H2}
it follows that for at least one $\eta$ for which
$\hat{U}^n\eta=\eta$ we have $\hat{T}^n\eta=\eta$. Hence by Lemma
~\ref{1.4} (ii) we have
$\eta=\hat{T}^n\eta=\eta_1^{kn}(\gamma_0)\eta$. Now using the facts
that the subgroup generated by $\gamma_0$ is dense in $G$ and that
$\hat{G}$ has no element of finite order we derive that
$\eta_1\equiv 1$. Hence $\eta\equiv 1$ which is a contradiction to
the fact that $\eta$ is non-trivial. So $(X,\mu,T)$ is ergodic.
\end{proof}

\begin{lemma}\label{1.6} If $\hat{G}$ has no element of
finite order, then for each $f\in C(X)$
\[\lim_{N\rightarrow\infty}\frac{1}{N}\sum_{n=0}^{N-1}f(T^n(x))\]
exists uniformly in $x$.
\end{lemma}
\begin{proof} Let $f\in C(X)$. Since $\hat{X}$ is dense in $C(X)$, assume, without loss of generality, that $f\in\hat{X}$.
Then $f=(f_1,f_2,\ldots,f_k)$ where $f_i\in\hat{G}$ for
$i=1,2,\ldots,k$. That is, for each $x=(x_1,x_2,\ldots,x_k)$ in $X$,
$f(x)=f_1(x_1)f_2(x_2)\cdots f_k(x_k)$. To examine the values
$f(T^n(x))=f(y_{1,n},y_{2,n},\ldots,y_{k,n})$, where
$y_{i,n}=(^{k}_{i})n^i\gamma_0+\sum_{j=1}^{i}(^{k-j}_{i-j})n^{i-j}x_j$
for $i=1,2,\ldots,k$, let
\[f_i(x_j)=\exp(\phi_i(x_j))\]
\[f_i(\gamma_0)=\exp(\alpha_i)\]
\noindent in which $0\leq\phi_i(x_j)<1$ and $0\leq\alpha_i<1$ for
$i,j=1,2,\ldots,k$. Then
 \[f(T^n(x))=f_1(y_{1,n})f_2(y_{2,n})\cdots
f_k(y_{k,n}).\] \noindent Since
\bea f_i(y_{i,n})&=&f_i((^{k}_{i})n^i\gamma_0+\sum_{j=1}^{i}(^{k-j}_{i-j})n^{i-j}x_j)\\
&=&\exp((^{k}_{i})n^i\alpha_i+\sum_{j=1}^{i}(^{k-j}_{i-j})n^{i-j}\phi_i(x_j))\eea

\noindent therefore we have
\bea f(T^n(x))&=&\exp(\sum_{i=1}^{k}(^{k}_{i})n^i\gamma_0+\sum_{i=1}^{k}\sum_{j=1}^{i}(^{k-j}_{i-j})n^{i-j}\phi_i(x_j))\\
&=&\exp(q(n))\eea \noindent where
\[q(n)=\sum_{i=1}^{k}(^{k}_{i})n^i\gamma_0+\sum_{i=1}^{k}\sum_{j=1}^{i}(^{k-j}_{i-j})n^{i-j}\phi_i(x_j).\]
Suppose that $f\neq 1$ and choose $j$ such that $1\leq j\leq k$,
$f_j\neq 1$ and $f_{j+1}=f_{j+2}=\ldots=f_{k}=1$. Then $q(n)$ is of
order $j$ and the leading coefficient in $q(n)$ equals to
$(^{k}_{j})\alpha_j$. Since $f_j$ is not of finite order and since
the subgroup generated by $\gamma_0$ is dense in $G$,
$\exp(\alpha_j)$ is not of finite order, that is, $\alpha_j$ is
irrational, so is $(^{k}_{j})\alpha_j$. Let
$q(n)=a_0+a_{1}n+\cdots+a_{j-1}n^{j-1}+(^{k}_{j})\alpha_j n^j$. Then
by Theorem ~\ref{Weyl},
\[\sum_{n=0}^{N-1}f(T^n(x))=\sum_{n=0}^{N-1}\exp(q(n))=\circ(N)\]
uniformly in $a_0,a_1,\ldots,a_{j-1}$. Since these $a_i$ are
combinations of $\phi_l(x_i)$, $\sum_{n=0}^{N-1}f(T^n(x))=\circ(N)$
uniformly in $x$. This proves the Lemma.
\end{proof}
Now we are ready to prove the main result of this section:

{\em Proof of Theorem ~\ref{1.7}.} (i) Since $(X,\mu,T)$ is ergodic
(Theorem ~\ref{1.5}), by the ergodic Theorem for each $f\in C(X)$
\[\lim_{N\rightarrow\infty}\frac{1}{N}\sum_{n=0}^{N-1}f(T^n(x))=\int_X
f d\mu\] for almost all $x$ in $X$. But it follows from Lemma
~\ref{1.6} that for
each $f\in C(X)$\\
\[\lim_{N\rightarrow\infty}\frac{1}{N}\sum_{n=0}^{N-1}f(T^n(x)):=f^\ast(x)\]
exists uniformly in $x$. So that $f^\ast$ is continuous. On the
other hand by the above equations $f^\ast$ is constant almost
everywhere, hence the continuity of $f^\ast$ implies that
\[\lim_{N\rightarrow\infty}\frac{1}{N}\sum_{n=0}^{N-1}f(T^n(x))=\int_X
f d\mu\] for all $x$ in $X$. Hence for each $x$ in $X$ the points
$T^n(x)$ ($n=1,2,3,\ldots$) are uniformly distributed (and so dense)
in $X$, that is, $(X,T)$ is minimal.

 (ii) Since $(X,T)$ is minimal
and $\hat{G}$ has no element of finite order, the result follows
from Lemma ~\ref{1.6} and Theorem~\ref{Oxtoby}.
\begin{remark}\label{r:1}
As a matter of fact, all the results of this section could be
generalized to arbitrary skew product dynamical systems. We avoided
doing this just because of the huge calculations needed.
\end{remark}

\section{Some consequences}
In this section, among other things, we deduce a simple proof for
the unique ergodicity of a class of functions on integers by using
the unique ergodicity of Milnes type skew product systems (Theorem
~\ref{1.7} (ii)). Also we show that Milnes type systems are natural
extensions of dynamical systems corresponding to some special
functions.\\
To this end, consider the homeomorphism $V:
l^{\infty}({\mathbb{Z}})\rightarrow l^{\infty}({\mathbb{Z}})$
defined by $V(f)(n)=f(n+1)$ (for $f\in l^{\infty}({\mathbb{Z}})$),
called the shift operator. For each $f\in l^{\infty}({\mathbb{Z}})$,
let $X_f$ be the pointwise closure of the set $\{V^n(f)$:
$n\in{\mathbb{Z}}\}$ in $l^{\infty}({\mathbb{Z}})$. Since on a
bounded subset of $l^{\infty}({\mathbb{Z}})$ the weak$^*$- topology
coincides with the topology of pointwise convergence on
${\mathbb{Z}}$, $X_f$ is compact with the pointwise topology. If we
denote the restriction of $V$ to $X_f$ by $V$ again, then the pair
$(X_f,V)$ defines a dynamical system. A function $f\in
l^{\infty}({\mathbb{Z}})$ is called uniquely ergodic if its
corresponding dynamical system $(X_f,V)$ is uniquely ergodic.

The following proposition was first proved in~\cite{S}, but for the
sake of completeness, we give the proof here.
\begin{proposition}\label{2.1} Let $(X,T)$ be a strictly ergodic
 dynamical system. Define $\Delta:C(X)\times X\rightarrow
 l^{\infty}({\mathbb{Z}})$ by $\Delta(F,x)(n)=F(T^n(x))$. Then any
 element in the range of $\Delta$ is uniquely ergodic.
 \end{proposition}
\begin{proof} Let $x_0\in X$ and $F\in C(X)$, and let
$f=\Delta(F,x_0)$. Define $\varphi:X\rightarrow
l^{\infty}({\mathbb{Z}})$ by $\varphi(x)(n)=F(T^n x)$. Then
$\varphi(x_0)=f$, and
$\varphi(Tx)(n)=F(T^n(Tx))=F(T^{n+1}x)=\varphi(x)(n+1)=V(\varphi(x))(n)$,
that is, $\varphi T=V\varphi$. Hence $\varphi$ is a flow
homomorphism. Since $(X,T)$ is minimal, $\varphi(X)=X_f$. It follows
that $f$ is uniquely ergodic, as claimed.
\end{proof}

Using Proposition~\ref{2.1} we derive a simple proof for the next
result which is due to Salehi~\cite{S} when $G={\mathbb{T}}$.

\begin{corollary}\label{2.3} Let $G$ be a compact monothetic group with
generator $\gamma_0$ such that $\hat{G}$ has no element of finite
order, then for every continuous character $\chi$ on $G$ and for
each $k\in{\mathbb{N}}$ the function $f(n)=\chi(n^k\gamma_0)$ is
uniquely ergodic. In particular for each $m$ in ${\mathbb{Z}}$ and
for each irrational element $\lambda$ in ${\mathbb{T}}$, the
function $f(n)=\lambda^{mn^k}$ is uniquely ergodic.
\end{corollary}
\begin{proof} Let $X=G^k$ and let $(X,T)$ be a Milnes type skew product dynamical system with the base group $G$, then by
Theorem~\ref{1.7} $(X,T)$ is strictly ergodic. Now
$T^n(0)=(kn\gamma_0,(^{k}_{2})n^2\gamma_0,\ldots,n^k\gamma_0)$. Let
$\chi\in\hat{G}$ and use Proposition ~\ref{2.1} with $x_0=0$ and the
function $F$ defined by $F(x_1,x_2,\ldots,x_k)=\chi(x_k)$. Then
$\Delta(F,0)=f$, where $f(n)=\chi(n^k\gamma_0)$. For the special
case, let $G={\mathbb{T}}$, $\chi=m$ and $\gamma_0=\lambda$, and use
the first part.
\end{proof}

Let $X={\mathbb{T}}^k$, let $\lambda\in{\mathbb{T}}$ and consider
the mapping $f:n\mapsto\lambda^{n^k}$ on ${\mathbb{Z}}$. Let the
mapping $\Gamma:X\rightarrow {\mathbb{T}}^{\mathbb{Z}}$ be defined
by
\[\Gamma(x_1,x_2,\ldots,x_k)(n)=f(n)x_1^{n^{k-1}}x_2^{n^{k-2}}\cdots
x_{k-1}^{n}x_k\] where $(x_1,x_2,\ldots,x_k)\in X$ and
$n\in{\mathbb{Z}}$. With these hypotheses we have the next result:

\begin{theorem}\label{2.4}
Let $(X,T)$ be a Milnes type skew product dynamical system with the
base group $G={\mathbb{T}}$. If $\lambda\in{\mathbb{T}}$ is
irrational, then $Range(\Gamma)=X_f$, and $\Gamma:(X,T)\rightarrow
(X_f,V)$ is a homomorphism, that is $V\Gamma=\Gamma T$ and so $X_f$
is a factor of ${\mathbb{T}}^k$.
\end{theorem}
\begin{proof}
Since $\lambda$ is irrational, by Theorem~\ref{1.7} $(X,T)$ is
minimal and hence the points
$(\lambda^{kn},\lambda^{(^{k}_{2})n^2},\ldots,\lambda^{n^k})=T^n(1,1,\ldots,1)$
($n=1,2,3,\ldots$) are dense in $X$. Now some computation reveals
that $Range(\Gamma)=X_f$ and also $V\Gamma=\Gamma T$.
\end{proof}

\begin{remarks}
(a) The mapping $\Gamma:(X,T)\rightarrow (X_f,V)$ in
Theorem~\ref{2.4} is not necessarily one-to-one, for example if
$k=3$, then $\Gamma(-1,-1,1)=\Gamma(1,1,1)$ because for each $n$,
$n^2+n$ is divisible by 2.

(b) Since any factor of any distal dynamical system is again distal,
it follows from Theorem~\ref{2.4} and Corollary~\ref{1.3} that for
each $k\in{\mathbb{N}}$ and for each $\lambda\in{\mathbb{T}}$ the
function $f:n\mapsto\lambda^{n^k}$ is a distal function, (note that
if $\lambda$ is not irrational, then $f$ is periodic and so distal).
This result is due to Namioka~\cite{N1}.
\end{remarks}

\section{Ellis Groups of skew product systems}

   In this section, unless otherwise stated, $(X,T)$ denotes a Milnes type
 skew product dynamical system with the base group $G$ and
generator $\gamma_0$ such that $\hat{G}$ has no element of finite
order, $1\neq k\in{\mathbb{N}}$, $\Sigma$ is the Ellis group of
$(X,T)$ and $E=E(G)$ (with its usual topology) is the (compact) set
of all (not necessarily continuous)
 endomorphisms of $G$.

 If $\tau=\lim_\alpha T^{m_\alpha}\in\Sigma$,
 then (passing to a subnet, if necessary) one may assume that for
 each $t$ in $G$ and for $i=1,2,\ldots,k-1$:
 \[\lim_\alpha T^{m_\alpha^i}(t)=\theta_i(t)\] and
\[\lim_\alpha T^{m_\alpha^k}\gamma_0=u\] If $x=(x_1,x_2,\ldots,x_k)\in
 X$ and $\tau x=(y_1,y_2,\ldots,y_k)$, then by Lemma \ref{1.2} (with $\theta_0=id_G$ in mind) for
 $i=1,2,\ldots,k-1$:
 \[y_i=\theta_i((^{k}_{i})\gamma_0)\prod_{j=1}^{i}\theta_{i-j}((^{k-j}_{i-j})x_j)\]
 and
\[y_k=u\prod_{j=1}^{k}\theta_{k-j}(x_j)\] Hence each
 $\tau$ in $\Sigma$ corresponds to a $k$-fold
 $(\theta_1,\theta_2,\ldots,\theta_{k-1},u)$ in $E^{k-1}\times
 G$. In the following Lemma it is shown that this correspondence is in
 fact an embedding isomorphism, in which the product in
 $E^{k-1}\times G$ is as
 follows:\\
 If $\tau\cong(\theta_{1},\ldots, \theta_{k-1}, u)$ and $\tau^{'}\cong({\theta}^{'}_{1},\ldots, {\theta}^{'}_{k-1},
u^{'})$, then $\tau^{'}\tau\cong(\phi_{1},\ldots, \phi_{k-1},z)$,
where for each $i=1,2,\ldots,k-1$,
$\phi_{i}=\prod_{j=0}^{i}{(^{i}_{j})}({\theta}^{'}_{i-j}o{\theta_{j}})$,
and
$z=u^{'}\prod_{j=1}^{k-1}{(^{k}_{j})}{\theta}^{'}_{k-j}\circ{\theta_{j}}(\gamma_0)u$.
For more information on these products see~\cite{M} and~\cite{N2}.

 The main theme of this section is to prove the following theorem which characterizes the topological
centre of $\Sigma$. Our method is partly similar to Example 3
in~\cite{M}. In what follows, for each $n$ in ${\mathbb{Z}}$, by
$n($ $)$ we mean the endomorphism $x\mapsto nx$ in $E(G)$.
\begin{theorem}\label{3.4} If $G$ is a connected
monothetic Lie group, then under the
identification given above:\\
\centerline {$\Lambda(\Sigma)=\{(n($ $),n^2($ $),\ldots,n^{k-1}($
$),u)\in\Sigma:n\in{\mathbb{Z}}$ and $u\in G\}\subset E^{k-1}\times
G$.}
\end{theorem}

To prove this theorem, we need to mention some preliminaries and
lemmas:

\begin{lemma}\label{3.1} Define $\Theta:\Sigma\rightarrow
E^{k-1}\times G$ by
$\Theta(\tau)=(\theta_1,\theta_2,\ldots,\theta_{k-1},u)$ where
$\tau\in\Sigma$ and $\theta_1,\theta_2,\ldots,\theta_{k-1},u$ are as
above. Then $\Theta$ is an embedding isomorphism of $\Sigma$ into
$E^{k-1}\times G$.
\end{lemma}
\begin{proof} $\Theta$ is well-defined. In fact, if for nets $\{m_{\alpha}\}$ and
$\{m_{\beta}^{'}\}$, $\tau=\lim_{\alpha}T^{m_{\alpha}}$,
$\tau^{'}=\lim_{\beta}T^{m_{\beta}^{'}}$ are members of $\Sigma$ and
$\tau=\tau^{'}$, we have to show that the associated elements
$\Theta(\tau)=(\theta_{1},\ldots, \theta_{k-1}, u)$ and
$\Theta(\tau^{'})=(\theta_{1}^{'},\ldots, \theta_{k-1}^{'}, u^{'})$
coincide in $E^{k-1}\times G$. If $x=(x_1,x_2,\ldots,x_k)\in X$,
then $\tau(x)=\tau^{'}(x)$, hence
\[u\prod_{j=1}^{k}\theta_{k-j}(x_j)=u^{'}\prod_{j=1}^{k}\theta^{'}_{k-j}(x_j)\]
By taking $x=0$, we obtain $u=u^{'}$. Now let $1\leq l\leq k-1$ and
let $(x_1,x_2,\ldots,x_k)\in X$ such that $x_l$ is an arbitrary
element of $G$ and for each $i\neq l$, $x_i=0$, the identity of $G$.
Then it follows that
$u\theta_{k-l}(x_l)=u^{'}\theta^{'}_{k-l}(x_l)$. But $u=u^{'}$ hence
$\theta_{k-l}=\theta^{'}_{k-l}$ and $\Theta$ is well-defined. That
$\Theta$ is a group isomorphism is clear from the formula. To prove
that $\Theta$ is continuous, let $\tau_\beta\rightarrow\tau$ in
$\Sigma$. Let $(x_1, x_2,..., x_k)\in X$ and let $\tau_\beta(x_1,
x_2,..., x_k)=(y_{1,\beta}, y_{2,\beta},..., y_{k,\beta})$ and
$\tau(x_1, x_2,..., x_k)=(y_1, y_2,..., y_k)$. If for each $\beta$,
$\Theta(\tau_\beta)=(\theta_{1,\beta},..., \theta_{(k-1),\beta},
u_\beta)$ and $\Theta(\tau)=(\theta_1,..., \theta_{k-1}, u)$, then
$y_{k,\beta}=u_{\beta}\prod_{j=1}^{k}\theta_{{k-j},\beta}(x_{j})$,
for all $\beta$, and $y_{k}=u\prod_{j=1}^{k}\theta_{k-j}(x_{j})$. If
we take $x_j=0$, for all $j=1, 2,...,k$, then since
$y_{k,\beta}\rightarrow y_k$, we get $u_\beta\rightarrow u$. Again
let $1\leq l\leq k-1$ and let $(x_1,x_2,\ldots,x_k)\in X$ such that
$x_l$ is an arbitrary element of $G$ and for each $i\neq l$,
$x_i=0$, then since $y_{k,\beta}\rightarrow y_k$ and
$u_\beta\rightarrow u$ we have $\theta_{k-l,\beta}\rightarrow
\theta_{k-l}$. Therefore $\Theta$ is an embedding isomorphism, as
claimed.
\end{proof}

\begin{lemma}\label{3.3} Let
$\tau\in\Lambda(\Sigma)$ and let $\Theta(\tau)=(\theta_1,\ldots,
\theta_{k-1}, u)$, where $\Theta$,
$\theta_1,\theta_2,\ldots,\theta_{k-1}$ and $u$ are as in Lemma
~\ref{3.1}, then for each $i=1,2,\ldots, k-1$, $\theta_i$ is
continuous.
\end{lemma}
\begin{proof} Let $\{t_\alpha\}$ be a net in $G$ and let
$t_\alpha\rightarrow t\in\ G$. Fix $1\leq i<k$. For each $\alpha$
let $x_\alpha$ (resp. $x$) be the element of $X$ in which all of its
coordinates are equal to $0$ except its $i$th coordinate which is
equal to $t_\alpha$ (resp. $t$). Then $x_\alpha\rightarrow x$ in
$X$. Let $\tau=\lim_\beta T^{m_\beta}\in\Lambda(\Sigma)$ and let
$\Theta(\tau)=(\theta_1,..., \theta_{k-1}, u)$. Assume that for each
$\alpha$, $\tau(x_\alpha)=(z_{1,\alpha}, z_{2,\alpha},...,
z_{k,\alpha})$ and $\tau(x)=(z_1, z_2,..., z_k)$. Since $(X,T)$ is
minimal (Theorem ~\ref{1.7} (i)) by Proposition ~\ref{m} $\tau$ is
continuous. Hence for each $j=1,2,\ldots,k$,
$z_{j,\alpha}\rightarrow z_j$, in particular,
$z_{k,\alpha}\rightarrow z_k$. But $z_{k,\alpha}=(\lim_\beta
T^{m_\beta^k}\gamma_0)\theta_{k-i}(t_\alpha)$ and $z_k=(\lim_\beta
T^{m_\beta^k}\gamma_0)\theta_{k-i} (t)$, hence
$\theta_{k-i}(t_\alpha)\rightarrow \theta_{k-i}(t)$, that is,
$\theta_{k-i}$ is continuous and the lemma is proved.
\end{proof}

{\em Proof of Theorem~\ref{3.4}.} Since $G$ is a connected
monothetic Lie group it is isomorphic (as a Lie group) to a torus
${\mathbb{T}}^m$ for some $m$ (see~\cite{DK}). Hence $\hat{G}$ has
no elements of finite order, therefore by Theorem ~\ref{1.7} (i)
$(X,T)$ is minimal. Now it follows from Proposition ~\ref{m} that
$\Lambda(\Sigma)=\Sigma_c$. To prove the theorem let
$\tau=\lim_\alpha T^{m_\alpha}$ be in $\Lambda(\Sigma)$ and let the
mapping $\Theta:\Sigma\rightarrow E^{k-1}\times G$ be defined as in
Lemma ~\ref{3.1} and let $\Theta(\tau)=(\theta_1,..., \theta_{k-1},
u)$. Hence by Lemma ~\ref{3.3} and the above results, $\theta_j$ is
continuous for each $j=1,2,\ldots,k-1$. Consider $\theta_j$ as an
endomorphism on ${\mathbb{T}}^m$. Then
$\theta_j(\eta_1,\eta_2,\ldots,\eta_m)=(f_j(\eta_1),f_j(\eta_2),\ldots,f_j(\eta_m))$
where $f_j(\eta)=\lim_\alpha\eta^{m_\alpha^j}$ for all $\eta$ in
${\mathbb{T}}$. The continuity of $\theta_j$ implies that $f_j\in
E({\mathbb{T}})$ is continuous, hence $f_j=($ $)^{n_j}$ for some
$n_j$ in ${\mathbb{Z}}$, and so $\theta_j=n_j($ $)$ on $G$. We will
show that for every $j=2,..., k-1$, $n_j($ $)=n_1^j($ $)$ on $G$.
Fix $p$ as a non-zero integer. Then for $\eta=\exp(\frac{1}{p})$,
$f_1(\eta)=\lim_\alpha(\eta)^{m_\alpha}= (\eta)^{n_1}$, hence
$m_\alpha=n_1$ (mod. $p$), eventually. Thus for each $j=2,..., k-1$,
$m_{\alpha}^j=n_1^j$ (mod. $p$), eventually. Therefore
$\lim_\alpha\exp(\frac{1}{p}m_{\alpha}^j)=\exp(\frac{1}{p}n_1^j)$.
It follows that for each integer $q$ with $0<q<p$ one has
$\lim_\alpha\exp(\frac{q}{p}m_{\alpha}^j)=\exp(\frac{q}{p}n_1^j)$.
On the other hand
$\lim_\alpha\exp(\frac{q}{p}m_{\alpha}^j)=\theta_j(\exp\frac{q}{p})$.
Thus ${\exp(\frac{q}{p})}^{n_1^j}={\exp(\frac{q}{p})}^{n_j}$. Since
$($ $)^{n_1^j}$ and $($ $)^{n_j}$ are both continuous functions on
${\mathbb{T}}$ and the set $\{\exp(\frac{q}{p}): p,q$ are positive
integers, $0<q<p\}$ is dense in ${\mathbb{T}}$ we have $($
$)^{n_1^j}=($ $)^{n_j}$. Thus $f_j=($ $)^{n_1^j}$ and so
$\theta_j=n_1^{j}($ $)$. The latter means that $\tau=(n_1($
$),n_1^2($ $),\ldots,n_1^{k-1},u)$. To prove the converse inclusion,
assume that $\tau=(n($ $ ), {n^2}($ $ ),..., {n^{k-1}}($ $
),u)\in\Sigma$ for some $n\in{\mathbb{Z}}$ and $u\in G$, then the
continuity of ${n^j}( $ $)$ for $j=1,2,\ldots,k-1$ implies that
$\tau$ is continuous, hence (by Proposition ~\ref{m})
$\tau\in\Lambda(\Sigma)$, and the theorem is proved.

\begin{remarks} \label{3.5}
(a) Let $(X,T)$ denote a Hahn's type skew product system. It would
be interesting to investigate the Ellis group $\Sigma$ and
$\Lambda(\Sigma)$ of the finite (i.e. $k<\infty$) or infinite (i.e.
$k=\infty$) skew product system $(X,T)$. Hahn~\cite{H1} proved that
if the base group $G$ has no element of finite order, then $(X,T)$
is a (distal) minimal dynamical system which is also uniquely
ergodic. It is readily proved by induction that for each $n$ in
${\mathbb{Z}}$ and for each $(x_1,x_2,\ldots,x_k)\in X$
\[{T}^n(x_1,x_2,\ldots,x_k)=(y_{1,n},y_{2,n},\ldots,y_{k,n})\]
where for each $j=1,2,\ldots k$
\[y_{j,n}=p_j(n)\gamma_0+p_{j-1}(n)x_1+p_{j-2}(n)x_2+\cdots+p_1(n)x_{j-1}+x_j\]
in which $p_j(n)=\frac{1}{j!}n(n-1)\cdots(n-j+1)$. The same formulas
work for the case in which $k=\infty$. Let $\Sigma$ be the Ellis
group corresponding to $(X,T)$. If $\tau\in\Sigma$ and for some net
$\{m_\alpha\}$ in ${\mathbb{Z}}$, $\tau =\lim_\alpha {T}^{m_\alpha}$
then by what we already discussed: if $k<\infty$, then $\tau$
corresponds to an element
$(\theta_1,\theta_2,\ldots,\theta_{k-1},u)$ in $E^{k-1}\times G$,
and if $k=\infty$, then $\tau$ corresponds to
$(\theta_1,\theta_2,\theta_3,\ldots)\in\prod_{j=1}^{\infty}E_j$ in
which for all $j$, $E_j=E(G)$ and for all $t\in G$,
$\theta_j(t)=\lim_\alpha p_j(m_\alpha)t$ (with $\theta_0=id_G$ in
mind), and $u=\lim_\alpha p_k(m_\alpha)\gamma_0$. It is not hard to
verify that under this correspondence if
$\tau\cong(\theta_{1},\ldots, \theta_{k-1}, u)$ and
$\tau^{'}\cong({\theta}^{'}_{1},\ldots, {\theta}^{'}_{k-1}, u^{'})$,
then $\tau^{'}\tau\cong(\phi_{1},\ldots, \phi_{k-1},z)$, where for
each $i=1,2,\ldots,k-1$,
$\phi_{i}=\prod_{j=0}^{i}{\theta}^{'}_{i-j}o{\theta_{j}}$ and
$z=u^{'}\prod_{j=1}^{k-1}{\theta}^{'}_{k-j}o{\theta_{j}}(\gamma_0)u$.
Similarly for the case $k=\infty$ if
$\tau\cong(\theta_{1},\theta_{2},\theta_{3},\ldots)$ and
$\tau^{'}\cong({\theta}^{'}_{1},{\theta}^{'}_{2},{\theta}^{'}_{3},\ldots)$,
then $\tau^{'}\tau\cong(\phi_{1},\phi_{2},\phi_{3},\ldots)$, where
for each $i$,
$\phi_{i}=\prod_{j=0}^{i}{\theta}^{'}_{i-j}o{\theta_{j}}$. With the
same methods as in Lemmas ~\ref{3.1} and ~\ref{3.3} one may easily
verify that if $k<\infty$, then the mapping
$\Theta:\Sigma\rightarrow E^{k-1}\times G$ defined by
$\Theta^{'}(\tau)=(\theta_1,\theta_2,\ldots,\theta_{k-1},u)$ where
$\theta_1,\theta_2,\ldots,\theta_{k-1}$ and $u$ are as above, is an
embedding isomorphism, and if $k=\infty$, then the mapping
$\Theta:\Sigma\rightarrow\prod_{j=1}^{\infty}E$,
$\tau\mapsto(\theta_1,\theta_2,\theta_3,\ldots)$ is an embedding
isomorphism. Also if $\tau\in\Lambda(\Sigma)$, then $\theta_i$ is
continuous for each $i=1,2,\ldots,k-1$, (and for all $i$, if
$k=\infty$).

 Now with the method of Theorem ~\ref{3.4} one
may observe that if $(X,T)$ is a finite or infinite Hahn's type skew
product dynamical system with the base group $G$ which is a
connected monothetic Lie group, then under the identifications given
above:

(i) If $k<\infty$, then\\
\centerline{$\Lambda(\Sigma)=\{(p_1(n)($ $),p_2(n)($
$),\ldots,p_{k-1}(n)($ $),u)\in\Sigma:n\in{\mathbb{Z}}$ and
$u\in G\}$}\\

(ii) If $k=\infty$, then $\Lambda(\Sigma)\equiv{\mathbb{Z}}$.\\

(b) In part (a) let k=2. Let $G={\mathbb{T}}$ and let
$\gamma_0\in{\mathbb{T}}$ be irrational. Namioka~\cite{N2} proved
that for this case $\Sigma\cong E({\mathbb{T}})\times{\mathbb{T}}$
and characterized the topological centre of $\Sigma$ explicitly. In
fact, he showed that
$\Lambda(\Sigma)\cong{\mathbb{Z}}\times{\mathbb{T}}=\{(($
$)^n,u):n\in{\mathbb{Z}}$ and $u\in{\mathbb{T}}\}$. So our results
generalize Namioka's example. Also, as discussed in~\cite{N2} for
the case $k=3$ in part (a), $\Sigma$ is not isomorphic to the whole
space $E({\mathbb{T}})\times E({\mathbb{T}})\times{\mathbb{T}}$. So
it would be interesting to determine which part of $E^{k-1}\times G$
is covered by $\Sigma$ (and also by $\Sigma$ in Theorem~\ref{3.4}).
There is a conjecture due to Milnes~\cite{M} concerning this problem
for the case when $G={\mathbb{T}}$ and $k=4$
in Theorem~\ref{3.4}.\\

(c) We do not know whether for an arbitrary
$k\in{\mathbb{N}}$,\\
 \centerline
{$\Lambda(\Sigma(X,T_M))\cong\{(n($ $),n^2($
$),\ldots,n^{k-1}($ $),u):n\in{\mathbb{Z}}, u\in G\}$ and}\\

\centerline {$\Lambda(\Sigma(X,T_H))\cong\{(p_1(n)($ $),p_2(n)($
$),\ldots,p_{k-1}(n)($ $),u):n\in{\mathbb{Z}}, u\in G\}$} \noindent
 where $(X,T_M)$ and $(X,T_H)$ denote any Milnes type and Hahn's type skew product dynamical system, respectively. To decide this using Theorem ~\ref{3.4} and
part (a), we have to characterize for which pairs $(n,u)$ in
${\mathbb{Z}}\times G$, one has $(n($ $),n^2($ $),\ldots,n^{k-1}($
$),u)\in\Sigma(X,T_M)$ and $(p_1(n)($ $),p_2(n)($
$),\ldots,p_{k-1}(n)($ $),u)\in\Sigma(X,T_H)$, respectively. Of
course (by part (b) and a similar method for $\Sigma(X,T_M)$) if
$k=2$ and $G={\mathbb{T}}$, then the answer to both questions is
affirmative. Hence another question
arises:\\
\noindent Is it true that if $k<\infty$, then
$\Lambda(\Sigma(X,T_M))\neq
{\mathbb{Z}}\neq\Lambda(\Sigma(X,T_H))$?\\
\noindent Note however that, as stated in part (a), for the case
$k=\infty$ we have $\Lambda(\Sigma(X,T_H))\equiv {\mathbb{Z}}$.\\

(d) Following similar methods for the Ellis groups of Hahn and
Milnes type systems, one may show that if $(X,T)$ is any skew
product dynamical system with the base group $G$ which is a
connected monothetic Lie group, then with our notations of section 3:\\

\centerline {$\Lambda(\Sigma(X,T))\cong\{(Q_1(n)($ $),Q_2(n)($
$),\ldots,Q_{k-1}(n)($ $),u)\in\Sigma(X,T):n\in{\mathbb{Z}}, u\in
G\}\cdot$}

 \noindent And also the product in $\Sigma(X,T)$ is given as
follows: If $\tau\cong(\theta_{1},\ldots, \theta_{k-1}, u)$ and
$\tau^{'}\cong({\theta}^{'}_{1},\ldots, {\theta}^{'}_{k-1}, u^{'})$,
then $\tau^{'}\tau\cong(\phi_{1},\ldots, \phi_{k-1},z)$, where for
each $i=1,2,\ldots,k-1$,
$\phi_{i}=\prod_{j=0}^{i}J(i,j){\theta}^{'}_{i-j}o{\theta_{j}}$ and
$z=u^{'}\prod_{j=1}^{k-1}J(k,j){\theta}^{'}_{k-j}o{\theta_{j}}(\gamma_0)u$.
(See also~\cite[Theorem 1]{M}.)\\

(e) The mapping $\Gamma:(X,T)\rightarrow (X_f,V)$ in
Theorem~\ref{2.4} induces a continuous homomorphism from $\Sigma$
onto the Ellis group $\Sigma({\mathbb{Z}},X_f)$ of $(X_f,V)$ (which
is not necessarily an isomorphism). We believe that it is worth
characterizing the topological centre of $\Sigma({\mathbb{Z}},X_f)$.
\end{remarks}

\section*{Acknowledgments} The first author would like to thank the Department of Mathematics and Statistics of the Alberta
University in Canada for their hospitality during his sabbatical
leaves in Alberta University, where this work was accomplished,
especially he would like to express his gratitude to Professor
Anthony Lau for his encouragement and support through his NSERC
grant A-7679. Also the authors would like to thank Professor Isaac
Namioka for his very many contributions to this work.

\end{document}